\documentclass[10pt,twoside,reqno]{amsart}
\usepackage{amssymb}
 \usepackage{amsmath}
\usepackage{multirow}
\usepackage{color}
\calclayout

\numberwithin{equation}{section}
\makeatletter

\renewcommand{\@secnumfont}{\bfseries}

\renewcommand{\section}{\@startsection{section}{1}%
  {0mm}{.7\linespacing\@plus\linespacing}{.5\linespacing}
  {\normalfont\bfseries\centering}}

\newcommand{\bibsection}{\@startsection{section}{1}%
  {0mm}{.7\linespacing\@plus\linespacing}{.5\linespacing}
  {\normalfont\scshape\centering}}

\renewcommand{\@biblabel}[1]{#1.}

\newtheorem{thm}{\bf Theorem}[section]

\newtheorem{cor}[thm]{\bf Corollary}

\begin{document}

\vspace{1.3cm}

\title{Degenerate zero-truncated Poisson random variables}

\author{Taekyun Kim}
\address{ Department of
Mathematics, Kwangwoon University, Seoul 139-701, Republic of
Korea}
\email{tkkim@kw.ac.kr}

\author{Dae San Kim}
\address{Department of Mathematics, Sogang University, Seoul 121-742, Korea}
\email{dskim@sogang.ac.kr}

\maketitle

\markboth{\centerline{\scriptsize Degenerate zero-truncated Poisson random variables}}{\centerline{\scriptsize T. Kim, D. S. Kim}}

\begin{abstract}
Recently, the degenerate Poisson random variable with parameter $\alpha > 0$, whose probability mass function is given by $P_{\lambda}(i) = e_{\lambda}^{-1} (\alpha)  \frac{\alpha^{i}}{i!} (1)_{i,\lambda}$, $(i = 0,1,2,\cdots)$, was studied.
In probability theory, the zero-truncated Poisson distributions are certain discrete probability distributions whose supports are the set of positive integers.
These distributions are also known as the conditional Poisson distributions or the positive Poisson distributions. In this paper, we introduce the degenerate zero-truncated Poisson random variables whose probability mass functions are a natural extension of the zero-truncated Poisson distributions, and investigate
various properties of those random variables.
\end{abstract}
\bigskip
\medskip

\section{\bf Introduction}
A random variable $X$ is a real valued function defined on a sample space. If $X$ takes any values in a countable set, then $X$ is called a discrete random variable.
For a discrete random variable $X$, the probability mass function $P(a)$ of $X$ is defined by
\begin{equation} \begin{split} \label{01}
P(a) = P\{X=a\},\,\,\,\,(\textnormal{see}\,\,[5,11,13,15,17]).
\end{split} \end{equation}

As is well known, a Poisson random variable indicates how many events occurred in a given period of time. A random variable $X$, taking on one of the
values $0,1,2,\cdots,$ is said to be a Poisson random variable with parameter $\alpha > 0$ if the probability mass function of $X$ is given by
\begin{equation} \begin{split} \label{02}
P(i) = P\{X=i\} = e^{-\alpha}  \frac{\alpha^{i}}{i!}, (i = 0,1,2,\cdots),\,\,\,\,(\textnormal{see}\,\,[17]).
\end{split} \end{equation}

The quantity $E[X^n]$ of the Poisson random variable with parameter $\alpha > 0$, which is
called the $n$-th moment of $X$, is defined by
\begin{equation} \begin{split} \label{03}
E[X^n] = \bigg(\sum_{i=0}^\infty i^n \frac{\alpha^i}{i!}\bigg) e^{-\alpha} = Bel_n(\alpha), (n \geq 1),\,\,\,\,
\end{split} \end{equation}
where $Bel_n(\alpha)$ are  the Bell polynomials defined by generating function as
\begin{equation} \begin{split} \label{04}
e^{\alpha (e^t-1)} = \sum_{n=0}^\infty Bel_n(\alpha) \frac{t^n}{n!} ,\,\,\,\,(\textnormal{see}\,\,[2,3,4,6,7,10,16]).
\end{split} \end{equation}

If $X$ is a discrete random variable taking values in the nonnegative integers, the probability
generating function of $X$ is defined as
\begin{equation}\label{05}
G(t) = E[t^X] = \sum_{n=0}^{\infty} P(n) t^n,
\end{equation}
where $P(n) = P\{X=n\}$ is the probability mass function of $X$.

Let $X = (X_1, X_2, \cdots, X_k)$ be a discrete random variable taking values in the $k$-dimensional
nonnegative integer lattice. Then the probability generating function of $X$ is given by
\begin{equation} \begin{split} \label{06-1}
G(t_1,t_2,\cdots,t_k) & = E[t_1^{X_1}, t_2^{X_2}, \cdots,t_k^{X_k}] \\
                      & = \sum_{x_1, x_2, \cdots, x_k = 0}^{\infty} P(x_1,x_2,\cdots,x_k) t_1^{x_1} t_2^{x_2} \cdots t_k^{x_k},
\end{split} \end{equation}
where $P(x_1,x_2,\cdots,x_k)$ is the probability mass function of $X$. The power series converges absolutely at least for all convex vectors $t = (t_1, t_2, \cdots,t_k) \in \mathbb{C}^k$ with
$max \{|t_1|, |t_2|, \cdots, |t_k| \} \leq 1$.

\vspace{0.1in}

For any real number $\lambda$, the $\lambda$-falling factorial sequence is defined by ({see, for more details, [1,8,9,18])
\begin{equation} \begin{split} \label{06}
(x)_{n,\lambda}=
\begin{cases}
  1, & \mbox{if } n=0, \\
  x(x-\lambda)(x-2\lambda)\cdots\big(x-(n-1)\lambda\big), & \mbox{if } n\geq 1.
\end{cases}
\end{split} \end{equation}

The degenerate exponential is defined by
\begin{equation} \begin{split} \label{07}
e_{\lambda}^{x}(t) = (1 + \lambda t)^{\frac{x}{\lambda}}, \quad e_{\lambda}(t) = e_{\lambda}^{1}(t),\,\,\,\,(\textnormal{see}\,\,[9,12]).
\end{split} \end{equation}

From \eqref{06} and \eqref{07}, we have
\begin{equation} \begin{split} \label{08}
e_{\lambda}^{x}(t) = \sum_{n=0}^{\infty}(x)_{n,\lambda} \frac{t^n}{n!},\,\,\,\,(\textnormal{see}\,\,[8,12]).
\end{split} \end{equation}

In [12], the degenerate Bell polynomials are given by the generating function 
\begin{equation} \begin{split} \label{09}
e_{\lambda}^{-1}(x) e_{\lambda}(xe^t) = \sum_{n=0}^\infty \beta_{n,\lambda}(x) \frac{t^n}{n!}.
\end{split} \end{equation}

From \eqref{04} and \eqref{09}, we note that
\begin{equation} \begin{split} \label{10}
\lim_{\lambda \rightarrow 0} \beta_{n,\lambda}(x) = Bel_n (x), \quad (n \geq 0).
\end{split} \end{equation}

For $\lambda \in \mathbb{R}$, $X_{\lambda}$ is the degenerate Poisson random variable with parameter $\alpha > 0$ if the
probability mass function of $X_{\lambda}$ is given by
\begin{equation} \begin{split} \label{11}
P_{\lambda} (i) = P \{X_{\lambda}= i \} = e_{\lambda}^{-1} {(\alpha)} \frac{\alpha^{i}}{i!} (1)_{i,\lambda},\,\,\,\,(\textnormal{see}\,\,[11]).
\end{split} \end{equation}
where $i = 0,1,2,\cdots$.

In probability theory, the zero-truncated Poisson distributions are certain discrete probability distributions whose supports are the set of positive integers. These distributions are also called the conditional Poisson distributions or the positive Poisson distributions.
A random variable $X$ is the zero-truncated Poisson random variable with parameter $\alpha >0$ if the probability mass function of $X$ is given by
\begin{equation} \begin{split} \label{12}
P(k) = P\{X=k\} =  \frac{ e^{-\alpha}}{1-e^{-\alpha}} \frac{\alpha^{k}}{k!} = \frac{1}{e^{\alpha}-1} \frac{\alpha^{k}}{k!},\,\,\,\,(\textnormal{see}\,\,[11,14]).
\end{split} \end{equation}
where $k = 1,2,3,\cdots$. 
In this paper, we introduce the degenerate zero-truncated Poisson random variables whose probability mass functions are a natural extension of the zero-truncated Poisson distributions, and investigate
various properties of those random variables.

\vspace{0.1in}
\section{Degenerate zero-truncated Poisson random variables}
\vspace{0.1in}

Let $X_{\lambda}$ be the discrete Poisson random variable with parameter $\alpha$.
Then we observe that
\begin{equation} \begin{split}
\sum_{n=1}^\infty P_{\lambda}(n) & = \sum_{n=1}^\infty \frac{\alpha^n}{n!} (1)_{n,\lambda}  e_{\lambda}^{-1} ({\alpha}) =   e_{\lambda}^{-1} ({\alpha}) \bigg(\sum_{n=0}^\infty \frac{\alpha^n}{n!} (1)_{n,\lambda} -1 \bigg) \\
& = e_{\lambda}^{-1}({\alpha}) \bigg(e_{\lambda}({\alpha}) -1 \bigg) = 1 - e_{\lambda}^{-1} ({\alpha}).
\end{split} \end{equation}

Thus, we note that
\begin{equation} \begin{split}\label{13}
\sum_{n=1}^\infty \frac{ P_{\lambda}(n)}{1 - e_{\lambda}^{-1} ({\alpha})} & = \sum_{n=1}^\infty \frac{\alpha^n (1)_{n,\lambda}} {1 - e_{\lambda}^{-1} ({\alpha})} \frac{e_{\lambda}^{-1} ({\alpha})}{n!} \\
& = \sum_{n=1}^\infty \frac{1}{e_{\lambda}({\alpha}) -1} \alpha^n \frac{(1)_{n,\lambda}}{n!} = 1.
\end{split} \end{equation}

From \eqref{13}, we can consider the degenerate zero-truncated Poisson random variable with parameter $\alpha > 0$.

A random variable $X_{\lambda}$ is the degenerate zero-truncated Poisson random variable with parameter $\alpha > 0$ if the probability mass function of $X$ is given by
\begin{equation} \begin{split} \label{14}
P_{\lambda}(n) = \frac{e_{\lambda}^{-1} (\alpha)}{1- e_{\lambda}^{-1} (\alpha)} \frac{\alpha^n}{n!} (1)_{n,\lambda} = \frac{1}{e_{\lambda}(\alpha)-1} \frac{\alpha^n}{n!} (1)_{n,\lambda},
\end{split} \end{equation}
where $n = 1,2,3,\cdots$.\\

Note that $\lim_{\lambda \rightarrow 0} P_{\lambda}(n) = \frac{1}{e^\alpha-1} \frac{\alpha^n}{n!} = P(n),\,\,\,(n \in \mathbb{N})$,
is the probability mass function of the zero-truncated Poisson random variable with parameter $\alpha > 0$.

The expectation of $X_{\lambda}$ is given by
\begin{equation} \begin{split}\label{15}
E[X_{\lambda}] & = \sum_{k=1}^\infty k P_{\lambda}(k) = \sum_{k=1}^\infty k \frac{\alpha^k}{k!} (1)_{k,\lambda} \frac{1}{e_{\lambda}(\alpha)-1} \\
               & = \frac{\alpha}{e_{\lambda}(\alpha)-1} \sum_{k=1}^\infty \frac{\alpha^{k-1}}{(k-1)!} (1)_{k-1,\lambda} (1- (k-1) \lambda)\\
               & = \frac{\alpha}{e_{\lambda}(\alpha)-1} \sum_{k=0}^\infty \frac{\alpha^{k}}{k!} (1)_{k,\lambda} - \frac{\alpha \lambda}{e_{\lambda}(\alpha)-1} \sum_{k=0}^\infty \frac{\alpha^{k}}{k!} k (1)_{k,\lambda} \\
               & = \frac{\alpha}{e_{\lambda}(\alpha)-1} e_{\lambda}(\alpha) - \alpha \lambda E[X_{\lambda}].
\end{split}\end{equation}

\noindent From \eqref{15}, we note that we have
\begin{equation}\begin{split}\label{16}
(1+ \alpha \lambda) E[X_{\lambda}] = \frac{\alpha}{e_{\lambda}(\alpha)-1} e_{\lambda}(\alpha).
\end{split} \end{equation}

\noindent Thus, by \eqref{16}, for $\lambda \neq -\frac{1}{\alpha}$ we get
\begin{equation} \begin{split} \label{17}
E[X_{\lambda}] = \frac{\alpha}{1+ \alpha \lambda} \frac{e_{\lambda}(\alpha)}{e_{\lambda}(\alpha)-1}.
\end{split} \end{equation}

\begin{thm}\label{thm1}
Let $X_{\lambda}$ be the degenerate zero-truncated Poisson random variable with parameter $\alpha > 0$.
Then, for $\lambda \neq -\frac{1}{\alpha}$, the expectation of $X_{\lambda}$ is given by
\begin{equation*}
\begin{split}
E[X_{\lambda}] & = \frac{\alpha}{1+ \alpha \lambda} \frac{ e_{\lambda}(\alpha)}{e_{\lambda}(\alpha)-1}\\
               & = \frac{\alpha}{1+ \alpha \lambda} \frac{1}{1- e_{\lambda}^{-1}(\alpha)}.
\end{split}
\end{equation*}
\end{thm}

Observing that $(1)_{k+2,\lambda}=(1)_{k+1,\lambda}-\lambda(k+1)(1)_{k+1,\lambda}$, we have
\begin{equation} \begin{split}\label{18}
E[X_{\lambda}^2] & = \sum_{k=1}^\infty k^2 P_{\lambda}(k) = \frac{1}{e_{\lambda}(\alpha)-1} \sum_{k=1}^\infty k^2 \frac{\alpha^k}{k!} (1)_{k,\lambda}  \\
               & = \frac{1}{e_{\lambda}(\alpha)-1} \sum_{k=2}^\infty \frac{\alpha^{k}}{(k-2)!} (1)_{k,\lambda} + \frac{1}{e_{\lambda}(\alpha)-1} \sum_{k=1}^\infty \frac{k}{k!} \alpha^{k} (1)_{k,\lambda} \\
               & = \frac{\alpha^2}{e_{\lambda}(\alpha)-1} \sum_{k=0}^\infty \frac{\alpha^{k}}{k!} (1)_{k+2,\lambda} + E[X_{\lambda}] \\
               & = \frac{\alpha^2}{e_{\lambda}(\alpha)-1} \left\{\sum_{k=0}^{\infty}\frac{\alpha^k}{k!}(1)_{k+1,\lambda}-\lambda\sum_{k=0}^{\infty}\frac{(k+1)\alpha^k}{k!}(1)_{k+1,\lambda}\right\}+ E[X_{\lambda}]\\
&=\frac{\alpha}{e_{\lambda}(\alpha)-1}\sum_{k=0}^{\infty}\frac{(k+1)\alpha^{k+1}}{(k+1)!}(1)_{k+1,\lambda}\\
&\quad\quad-\frac{\alpha\lambda}{e_{\lambda}(\alpha)-1}\sum_{k=0}^{\infty}\frac{(k+1)^2\alpha^{k+1}}{(k+1)!}(1)_{k+1,\lambda}+E[X_{\lambda}].\\         
\end{split}\end{equation}

Thus, by \eqref{18},  we get
\begin{equation}\label{19}
E[X_{\lambda}^2]= \alpha E[X_{\lambda}]-\alpha\lambda E[X_{\lambda}^2]+E[X_{\lambda}].
\end{equation}

and hence, from \eqref{19} and Theorem 1, for $\lambda \neq -\frac{1}{\alpha}$ we have
\begin{equation} \label{20}
E[X_{\lambda}^2] = \frac{1+\alpha}{1+\alpha\lambda}E[X_{\lambda}]=\frac{\alpha(1+\alpha)}{(1+\alpha\lambda)^2}\frac{e_{\lambda}(\alpha)}{e_{\lambda}(\alpha)-1}.
\end{equation}

For $\lambda \neq -\frac{1}{\alpha}$, the variance of $X_{\lambda}$ is given by
\begin{equation} \begin{split}\label{21}
Var[X_{\lambda}] = &E[(X_{\lambda}-E[X_{\lambda}])^2] = E[X_{\lambda}^2] - (E[X_{\lambda}])^2 \\
=& \frac{\alpha(1+\alpha)}{(1+\alpha\lambda)^2}\frac{e_{\lambda}(\alpha)}{e_{\lambda}(\alpha)-1}
-\Big(\frac{\alpha}{1+\alpha\lambda}\Big)^2 \Big(\frac{e_{\lambda}(\alpha)}{e_{\lambda}(\alpha)-1}\Big)^2.
\end{split}\end{equation}

Therefore, by \eqref{21}, we obtain the following theorem.

\begin{thm}\label{thm2}
Let $X_{\lambda}$ be the degenerate zero-truncated Poisson random variable with parameter $\alpha > 0$.
Then, for $\lambda \neq -\frac{1}{\alpha}$, we have
\begin{equation*}
Var[X_{\lambda}]=\frac{\alpha(1+\alpha)}{(1+\alpha\lambda)^2}\frac{e_{\lambda}(\alpha)}{e_{\lambda}(\alpha)-1}
-\Big(\frac{\alpha}{1+\alpha\lambda}\Big)^2 \Big(\frac{e_{\lambda}(\alpha)}{e_{\lambda}(\alpha)-1}\Big)^2.
\end{equation*}
\end{thm}

Note that
\begin{equation*} \begin{split}
\lim_{\lambda \rightarrow 0} Var[X_{\lambda}] & = \alpha(1+\alpha)\frac{e^{\alpha}}{e^{\alpha}-1}-\alpha^2\Big(\frac{e^{\alpha}}{e^{\alpha}-1}\Big)^2 \\
& = \frac{\alpha e^{\alpha}}{(e^{\alpha}-1)^2} \big(e^{\alpha}-\alpha-1 \big).
\end{split}\end{equation*}
is the variance of the zero-truncated Poisson random variable with parameter $\alpha > 0$.

Now, we consider the generating function of the moment of the degenerate zero-truncated Poisson random variable with $\alpha > 0$.
From \eqref{14}, we note that
\begin{equation} \begin{split}\label{22}
E[e^{X_{\lambda} t}] & = \sum_{n=1}^\infty e^{nt} P_{\lambda}(n) = \frac{1}{e_{\lambda}(\alpha)-1} \sum_{n=1}^\infty \frac{\alpha^n}{n!} (1)_{n,\lambda} e^{nt}  \\
               & = \frac{1}{e_{\lambda}(\alpha)-1} \bigg\{ \sum_{n=0}^\infty \frac{\alpha^n}{n!} (1)_{n,\lambda} e^{nt} -1 \bigg\} \\
               & = \frac{1}{e_{\lambda}(\alpha)-1} \bigg(e_{\lambda}(\alpha e^t) -1 \bigg).
\end{split}\end{equation}

By \eqref{09} and \eqref{22}, we get
\begin{equation} \begin{split}\label{23}
\sum_{n=0}^\infty E[X_{\lambda}^n] \frac{t^n}{n!} & = E[e^{X_{\lambda} t}] = \frac{1}{e_{\lambda}(\alpha)-1} \bigg(e_{\lambda}(\alpha e^t) -1 \bigg) \\
 & = \frac{e_{\lambda}(\alpha)}{e_{\lambda}(\alpha)-1} \bigg(e_{\lambda}^{-1}(\alpha) e_{\lambda}(\alpha e^t) - e_{\lambda}^{-1}(\alpha) \bigg) \\
& = \frac{1}{1- e_{\lambda}^{-1}(\alpha)} \bigg(\sum_{n=0}^\infty \beta_{n,\lambda}(\alpha) \frac{t^n}{n!} - e_{\lambda}^{-1}(\alpha) \bigg) \\
& = \frac{1}{1- e_{\lambda}^{-1}(\alpha)} \sum_{n=1}^\infty \beta_{n,\lambda}(\alpha) \frac{t^n}{n!} + \frac{1}{1- e_{\lambda}^{-1}(\alpha)} (1- e_{\lambda}^{-1}(\alpha)) \\
  & = 1 + \frac{1}{1- e_{\lambda}^{-1}(\alpha)} \sum_{n=1}^\infty \beta_{n,\lambda}(\alpha) \frac{t^n}{n!}.
\end{split}\end{equation}

Therefore, by comparing the coefficients on both sides of \eqref{23}, we obtain the following theorem.

\begin{thm}\label{thm3}
Let $X_{\lambda}$ be the degenerate zero-truncated Poisson random variable with parameter $\alpha > 0$.
Then we have
\begin{equation*}
\begin{split}
E[X_{\lambda}^n] = \frac{1}{1- e_{\lambda}^{-1}(\alpha)} \beta_{n,\lambda}(\alpha),\,\,\,\,(n \in \mathbb{N}).
\end{split}
\end{equation*}
In particular,
\begin{equation*}
\begin{split}
E[1] = E[X_{\lambda}^0] = 1.
\end{split}
\end{equation*}
\end{thm}

We note that, for $x \geq 1$, the cumulative distribution function is given by
\begin{equation} \begin{split}\label{24}
F_{X_{\lambda}}(x) & = P \{X_{\lambda} \leq x \} =  \sum_{k=1}^{[x]} P_{\lambda}(k) \\
  & = \frac{1}{e_{\lambda}(\alpha)-1} \sum_{k=1}^{[x]} \frac{\alpha^k}{k!} (1)_{k,\lambda} \\
   & = \frac{1}{e_{\lambda}(\alpha)-1} \bigg(\sum_{k=0}^{[x]} \frac{\alpha^k}{k!} (1)_{k,\lambda} -1 \bigg),
\end{split}\end{equation}
where $[x]$ is the integral part of $x$.

Let us define the function $e_{\lambda,b}(a)$ to be
\begin{equation} \begin{split}\label{25}
e_{\lambda,b}(a) = \sum_{k=0}^b \frac{a^k}{k!} (1)_{k,\lambda},\,\,\,\,(b \in \mathbb{N}).
\end{split}\end{equation}

Therefore, by \eqref{24} and \eqref{25}, we obtain the following theorem.

\begin{thm}\label{thm4}
Let $X_{\lambda}$ be the degenerate zero-truncated Poisson random variable with parameter $\alpha > 0$.
Then the cumulative distribution function of $X_{\lambda}$ is given by
\begin{equation*}
\begin{split}
F_{X_{\lambda}}(x) = \frac{1}{e_{\lambda}(\alpha)-1} \bigg(e_{\lambda,[x]}(\alpha) -1 \bigg).
\end{split}
\end{equation*}
\end{thm}

As is known, the Stirling number of second kind is defined by
\begin{equation} \begin{split} \label{26}
x^n = \sum_{k=0}^n S_2(n,k) (x)_k, \quad (n \geq 0),\quad(\textnormal{see}\,\,[6,14,16,18]),
\end{split} \end{equation}
where $(x)_0=1,\,\, (x)_k = x(x-1) \cdots (x-k+1),\,\,(k \geq 1).$

From \eqref{26}, we can easily get
\begin{equation} \begin{split} \label{27}
\frac{1}{k!} (e^t-1)^k = \sum_{n=k}^\infty S_2(n,k) \frac{t^n}{n!},\,\,\,\,(k \geq 0).
\end{split} \end{equation}

In [9], the degenerate Stirling number of the second kind is defined by the generating function 
\begin{equation} \begin{split} \label{28}
\frac{1}{k!} (e_{\lambda}(t)-1)^k & = \frac{1}{k!} \bigg((1+ \lambda t)^{\frac{1}{\lambda}} -1 \bigg)^k \\
                                  & = \sum_{n=k}^\infty S_2(n,k) \frac{t^n}{n!},\,\,\,\,(k \geq 0).
\end{split} \end{equation}

From \eqref{28}, we note that
\[
S_{2,\lambda}(n+1,k)=kS_{2,\lambda}(n,k)+S_{2,\lambda}(n,k-1)-n\lambda S_{2,\lambda}(n,k),
\]
where $1\leq k \leq n$.

Suppose that $X_{1,\lambda}, X_{2,\lambda},\cdots, X_{k,\lambda}$ are independent degenerate zero-truncated Poisson random variables with parameter $\alpha(>0)$.
Let
\begin{equation}\label{29}
  X_\lambda=\sum_{i=1}^k X_{i,\lambda}, ~(k\in\mathbb{N}).
\end{equation}
From \eqref{05} and \eqref{11}, we note that
\begin{equation}\label{30}
  \begin{aligned}
  E[t^{X_{i,\lambda}}]=&\sum_{x=1}^\infty P_\lambda[X_{i,\lambda}= x]t^x\\
  =&\frac{1}{e_\lambda(\alpha)-1}\sum_{n=1}^\infty \frac{\alpha^n}{n!} (1)_{n,\lambda} t^n\\
  =&\frac{1}{e_\lambda(\alpha)-1}\left( \sum_{n=0}^\infty \frac{\alpha^n}{n!} (1)_{n,\lambda} t^n-1\right) \\
  =&\frac{e_\lambda(\alpha t)-1}{e_\lambda(\alpha)-1}.
  \end{aligned}
\end{equation}

By \eqref{29}, we get
\begin{equation}\label{31}
\begin{aligned}
  E[t^{X_{\lambda}}]=& E[t^{X_{1,\lambda}+X_{2,\lambda}+\cdots+X_{k,\lambda}}]\\
  =&\prod_{i=1}^k E[t^{X_{i,\lambda}}]\\
  =&\left( \frac{e_\lambda(\alpha t)-1}{e_\lambda(\alpha)-1}\right)^k.
\end{aligned}
\end{equation}
From \eqref{28} and \eqref{31}, we can derive the following equation:
\begin{equation}\label{32}
\begin{aligned}
  E[t^{X_{\lambda}}]=& \left(\frac{1}{e_\lambda(t)-1}\right)^k k! \frac{1}{k!}\left({e_\lambda(\alpha t)-1}\right)^k\\
  =&\frac{k!}{(e_\lambda(\alpha)-1)^k} \sum_{n=k}^\infty S_{2,\lambda}(n,k) \alpha^n \frac{t^n}{n!}\\
  =&\sum_{n=k}^\infty \left\{ \frac{k!}{(e_\lambda(\alpha)-1)^k}\frac{\alpha^n}{n!} S_{2,\lambda}(n,k)\right\} t^n.
\end{aligned}
\end{equation}

On the other hand, by \eqref{05}, we get
\begin{equation}\label{33}
\begin{aligned}
  E[t^{X_{\lambda}}]=& E\big[t^{X_{1,\lambda}+X_{2,\lambda}+\cdots+X_{k,\lambda}}\big]\\
  =&\sum_{n=k}^\infty P_\lambda\Big[\sum_{i=1}^k X_{i,\lambda}=n\Big]t^n.
\end{aligned}
\end{equation}
Therefore, by \eqref{32} and \eqref{33}, we obtain the following assertion.

\begin{thm}
  Let $X_{1,\lambda}, X_{2,\lambda}, \cdots, X_{k,\lambda}$ be independent degenerate zero-truncated Poisson random variables with parameter $\alpha(>0)$, and let $X_\lambda=\sum\limits_{i=1}^k X_{i,\lambda}$. Then the probability for $X_\lambda$ is given by
\begin{equation*}
P_\lambda[X_\lambda=n]=\begin{cases}
\frac{k!}{(e_\lambda(\alpha)-1)^k}\frac{\alpha^n}{n!}S_{2,\lambda}(n,k), & \mbox{if } n\geq k, \\
 0, & \mbox{otherwise}.
\end{cases}
\end{equation*}
\end{thm}

By \eqref{07}, we easily get
\begin{equation}\label{34}
\begin{aligned}
\frac{1}{k!} (e_\lambda(t)-1)^k=& \frac{1}{k!} \sum_{l=0}^k \binom{k}{l} (-1)^{k-l} e_\lambda^l (t)\\
=&\frac{1}{k!} \sum_{l=0}^k \binom{k}{l} (-1)^{k-l} \sum_{n=0}^\infty (l)_{n,\lambda} \frac{t^n}{n!}\\
=&\sum_{n=0}^\infty \left(\frac{1}{k!} \sum_{l=0}^k \binom{k}{l}(-1)^{k-l} (l)_{n,\lambda}\right) \frac{t^n}{n!}.
\end{aligned}
\end{equation}

From \eqref{28} and \eqref{34}, we have
\begin{equation}\label{35}
  \frac{1}{k!} \sum_{l=0}^k \binom{k}{l} (-1)^{k-l}(l)_{n,\lambda}=
  \begin{cases}
    S_{2,\lambda}(n,k), & \mbox{if } n\geq k, \\
    0, & \mbox{otherwise}.
  \end{cases}
  \end{equation}
Therefore, by \eqref{35} and Theorem 2.5, we obtain the following corollary.

\begin{cor}
Let $X_{1,\lambda}, X_{2,\lambda}, \cdots, X_{k,\lambda}$ be independent degenerate zero-truncated Poisson random variables with parameter $\alpha(>0)$ and let $X_\lambda=\sum\limits_{i=1}^k X_{i,\lambda}$. Then the probability for $X_\lambda$ is given by
\begin{equation*}
\begin{aligned}
P_\lambda[X_\lambda=n]=&\frac{\alpha^n}{n!(e_\lambda (\alpha)-1)^k}\sum_{l=0}^k \binom{k}{l} (-1)^{k-l} (l)_{n,\lambda}\\
=&\frac{\alpha^n}{n!(e_\lambda (\alpha)-1)^k}\sum_{l=0}^k \binom{k}{l} (-1)^{l} (k-l)_{n,\lambda},
\end{aligned}
\end{equation*}
where $n,~k\in\mathbb{N}$ with $n\geq k$.
\end{cor}

Assume that $X_{1,\lambda}, X_{2,\lambda}, \cdots, X_{k,\lambda}$ are independent degenerate zero-truncated Poisson random variables with respective parameters $\alpha_1, \alpha_2, \cdots, \alpha_k$. We let $X_\lambda=\sum\limits_{i=1}^k X_{i,\lambda}$. Then, by \eqref{05}, we get
\begin{equation}\label{36}
\begin{aligned}
\sum_{n=k}^\infty P_\lambda[X_\lambda=n]t^n=&E[t^{X_\lambda}]\\
=&E\left[ t^{X_{1,\lambda}+\cdots+X_{k,\lambda}}\right]\\
=&\prod_{i=1}^k E[t^{X_{i,\lambda}}].
\end{aligned}
\end{equation}

From \eqref{11}, we have
\begin{equation}\label{37}
\begin{aligned}
 E[t^{X_{i,\lambda}}]=&\sum_{x=1}^\infty P[X_{i,\lambda}=x]t^x\\
 =&\frac{1}{e_\lambda(\alpha_i)-1} \sum_{n=1}^\infty \frac{\alpha_i^n}{n!} (1)_{n,\lambda}t^n.
\end{aligned}
\end{equation}

By \eqref{36} and \eqref{37}, we obtain
\begin{equation}\label{38}
\begin{aligned}
 E[t^{X_{\lambda}}]=&\prod_{i=1}^k E[t^{X_{i,\lambda}}]\\
 =&\left( \sum_{n_1=1}^\infty \frac{\alpha_1^{n_1}}{n_1!}(1)_{n_1,\lambda}t^{n_1}\right)\cdots
 \left( \sum_{n_k=1}^\infty \frac{\alpha_k^{n_k}}{n_k!}(1)_{n_k,\lambda}t^{n_k}\right)\prod_{i=1}^k
 \frac{1}{e_\lambda(\alpha_k)-1}\\
  =&\sum_{n=k}^\infty\prod_{i=1}^k\left( {e_\lambda(\alpha_i)-1}\right)^{-1}\sum_{n_1+\cdots+n_k=n}\binom{n}{n_1,\cdots,n_k}
   \left({\alpha_1^{n_1}}(1)_{n_1,\lambda}\right)\cdots\\
   &\times
 \left({\alpha_k^{n_k}}(1)_{n_k,\lambda}\right) \frac{t^n}{n!}.
\end{aligned}
\end{equation}

Therefore, by \eqref{36} and \eqref{38}, we get
\begin{align*}
P_\lambda[X_\lambda=n]=&\frac{1}{n!}
\left( \prod_{i=1}^k \frac{1}{e_\lambda(\alpha_i)-1}\right)\sum_{n_1+\cdots+n_k=n}\binom{n}{n_1,\cdots,n_k}
   \left({\alpha_1^{n_1}}(1)_{n_1,\lambda}\right)\cdots\\
   &\times
 \left({\alpha_k^{n_k}}(1)_{n_k,\lambda}\right),
\end{align*}
where $n,~k\in \mathbb{N}$ with $n\geq k$.

\section{\bf Conclusions}
It has been shown that various probabilistic methods can be applied to the study of some special numbers and polynomials arising from combinatorics and number theory. For example, see [11,13].\\
\indent Let $X_{\lambda}$ be the degenerate zero-truncated Poisson random variable with parameter $\alpha$.
Then, for the random variable $X_{\lambda}$, we derived its expectation, its variance, its $n$-th moment, and its cumulative distribution function. In addition, we obtained two different expressions for the probability function of a finite sum of independent degenerate zero-truncated Poisson random variables with equal parameters. Furthermore, we were able to get an expression for the probability function of a finite sum of independent degenerate zero-truncated Poisson random variables with unequal parameters. \\
\indent As one of our future projects, we would like to continue this line of research, namely to explore applications of various methods of probability theory to the study of some special polynomials and numbers.

\end{document}